\documentclass[12pt,a4paper]{article}
\usepackage{amssymb,amsthm,amsmath}
\usepackage{fullpage}

\newtheorem{e-proposition}[theorem]{Proposition}

\newtheorem{e-definition}[theorem]{Definition\rm}

\newtheorem{theoreme}{Th\'eor\`eme}[section]

\newtheorem{proposition}[theoreme]{Proposition}

\DeclareMathSymbol\crossrt{\mathrel}{AMSb}{"6E}
\def\N{\mathbb N}
\def\R{\mathbb R}

\def\CU{{\mathcal U}}
\def\CA{{\mathcal A}}
\def\CH{{\mathcal H}}
\def\CK{{\mathcal K}}

\def\oh{\frac{1}{2}}

\def\oq{\frac{1}{q}}
\def\acts{\triangleright}
\def\ket#1{|#1\rangle}
\def\sq2{\CA(S^2_q)}
\def\suq2{\CU_q(su(2))}
\newcommand{\nn}{\nonumber}
\begin{document}

\title{The Spectral Geometry \\ of the \\ Equatorial Podle\'s Sphere} 
\date{}
\author{~\\ Ludwik D\c{a}browski $\strut^{1}$, 
Giovanni Landi $\strut^{2}$, \\ [15pt] 
Mario Paschke $\strut^{3}$, 
Andrzej Sitarz $\strut^{4}$ \thanks{Partially supported by Polish State
Committee for Scientific Research (KBN) under grant 2 P03B 022 25.} \\[30pt]
$\strut^{1}$ {\small Scuola Internazionale Superiore di Studi Avanzati,}\\[5pt]
{\small Via Beirut 2-4, I-34014, Trieste, Italy} \\[8pt] 
$\strut^{2}$ {\small Dipartimento Matematica e Informatica, Universita di Trieste}
\\[5pt] {\small via A. Valerio 12/1, I-34127, Trieste, Italy} \\[8pt]
$\strut^{3}$ {\small Max-Planck-Institut f\"ur Mathematik in den
Naturwissenschaften,}\\[5pt] {\small Inselstr.\ 22, 04103 Leipzig, Germany} \\[8pt] 
$\strut^{4}$ {\small Institute of Physics, Jagiellonian University,} \\[5pt] {\small
Reymonta 4, 30-059 Krak\'ow, Poland} }
\maketitle

\begin{abstract}
We propose a slight modification of the properties
of a spectral geometry a la Connes, 
allowing for some of the algebraic relations to be satisfied modulo
compact operators. On the equatorial Podle\'s sphere  we construct
$\suq2$-equivariant Dirac operator and real structure which
satisfy these modified properties.  
\end{abstract}

\vspace{2pc}



\vfill\eject

\section{Introduction}
Recent examples of noncommutative spectral geometries on spaces
coming from quantum groups \cite{ChPa,DaSi,Co,SW,NT,Kr} 
opened a new interesting and promising research area.
Along these lines, introducing a slight modification of the
defining properties of a noncommutative geometry, we present
a spectral geometry for the equatorial
Podle\'s sphere $S_q^2$ of \cite{Po}. Both the Dirac operator $D$ and the
real structure $J$ are equivariant under the action of
$\suq2$ on $S_q^2$.

With $q$ a real number, $0<q \leq 1$, we denote by $\sq2$ the algebra of polynomial
functions generated by operators $a$, $a^*$ and $b=b^*$, which
satisfy the following commutation rules: 
$$
ba = q^2 ab, \qquad\; a^*b = q^2 ba^*, \qquad\; a^* a + b^2 =1, \qquad\; 
q^2 a a^* + q^{-2} b^2 = q^2.
$$ 
This algebra has an  `equator'-worth of classical points
given by the one dimensional representations
$b=0, a=\lambda$, with $\lambda \in S^1$.

The symmetry of the sphere, which we shall use for the equivariance, is the
Hopf algebra module structure with respect to the $\suq2$ Hopf
algebra. It is derived from the canonical $\suq2$ action on the
$\CA(SU_q(2))$ algebra. Explicitly, the generators $e,f,k$
of $\suq2$ act on the generators of $\sq2$ in the following
way:
\begin{align} &k \acts a = q a, \;\;\;\; & &e \acts a = - (1+q^2) q^{-\frac{5}{2}} b,
\;\;\;& &f
\acts a = 0, \nn \\ &k \acts a^* = \oq {a^*}, & &e \acts {a^*} = 0, & &f \acts {a^*}
= (1+q^2) q^{-\frac{3}{2}} b, \label{module} \\ &k \acts b = b,	& &e \acts b =
q^{\frac{1}{2}} {a^*}, & &f
\acts b = -q^{\frac{3}{2}} a. \nn
\end{align}

The irreducible finite dimensional
representations of the Hopf algebra $\suq2$ are labelled by
a positive half-integers (see \cite{Maj}, e.g.) and each
representation space $V_l$ has a basis
$\{\ket{l,m}, \, m\in\{-l,-l+1,\dots ,l\}$ declared
to be orthonormal.
\section{Variations on spectral geometry}
A spectral geometry (a la Connes) are data
$(\CA,\pi,\CH,\gamma,J,D)$ fulfilling a series of
requirements \cite{grav}.

On the equatorial sphere $\CA=\CA(S^2_q)$
we construct an equivariant spectral geometry \cite{AS1},
starting from an equivariant representation $\pi$ on a
suitable Hilbert space $\CH$. On the latter there are
an equivariant real structure $J$ and an equivariant
Dirac operator $D$. However, with such a $J$ it is not possible
to satisfy all the requirements of \cite{grav}. Nevertheless, the algebraic requirements shall be obeyed up
to compact operators. In particular, the antilinear isometry
$J$, which provides the real structure, will map $\pi(\CA)$
to its commutant only modulo compact operators.
\begin{equation}
\label{alcomm}
\forall x,y\in \CA: \;\;\; \left[ \pi(x), J\pi(y)J^{-1} \right] \in  \CK \,,
\end{equation}
and the first order condition is valid only modulo compact operators,
\begin{equation}
\label{alFOC}
\forall x,y\in \CA: \;\;\; \left[ J\pi(x)J^{-1}, [D, \pi(y)] \right]
\in \CK \, .
\end{equation}
The essentially unique Dirac operator, which comes out as the solution
of the above condition, shall have the crucial property that all
commutators $[D,\pi(x)]$, $x \in \sq2$, are bounded.
\section{The equivariant geometry of $\sq2$}
A fully equivariant approach, both for the real structure
and the Dirac operator was worked out for the standard Podle\'s
sphere in \cite{DaSi}. Here we build up on this and an earlier
approach \cite{LPS}  
(see also \cite{MPTh}) to construct an equivariant spectral
geometry of the Equatorial Podle\'s sphere.

\bigskip
The starting ingredient of the equivariant spectral geometry, which
we are about to construct is an equivariant representation of
$\sq2$ on a Hilbert space $\CH$. Let us recall that an $H$-equivariant
representation of an $H$-module algebra $\CA$ on $\CH$, is
a representation $\pi$ of $\CA$ on $\CH$ and a representation
$\rho$ of $H$ on a dense subspace of $\CH$ such that for every
$a \in \CA$ and $h \in H$ and $v$ from a dense subspace we have:
\begin{equation}
\rho(h) \pi(a) v = \pi(h_{(1)} \triangleright a) \rho (h_{(2)}) v. 
\label{equiv-1}
\end{equation}
In our case, this is the same as a representation of the crossed
product of ~\mbox{$\suq2\crossrt\sq2$} on a dense subspace of the
Hilbert space. The general theory for the entire family of
Podle\'s spheres is in \cite{SW1}.
\begin{proposition}\label{piequi}
There exists two irreducible $\suq2$-equivariant representations
(denoted by $\pi_{\pm}$) of the algebra $\sq2$ on the Hilbert
space
$\CH_h = \bigoplus_{l=\frac{1}{2},\frac{3}{2},\ldots} V_l$
given by :
\begin{align}
\pi_\pm(a) \ket{l,m} =& \pm ( 1+q^2) \frac{q^{m-\oh}}{[2l][2l+2]}
~\sqrt{[l+m+1][l-m]}\; \ket{l,m+1}
\nn \\ &+ \frac{q^{m-l-\oh}}{[2l+2]} ~\sqrt{ [l+m+1][l+m+2]} \;\ket{l+1,m+1} \nn \\ &-
\frac{q^{m+l+\oh}}{[2l]} ~\sqrt{ [l-m][l-m-1]} \;\ket{l-1,m}, \nn \\ ~ \nn \\
\pi_\pm(b) \ket{l,m} =& \pm \frac{1}{[2l][2l+2]} \left( [l-m+1][l+m] - q^2
[l-m][l+m+1] \right)\;
\ket{l,m} \nn \\ & - \frac{q^{m+1}}{[2l+2]} ~\sqrt{ [l-m+1][l+m+1]}\; \ket{l+1,m} \nn
\\ &-
\frac{q^{m+1}}{[2l]} ~\sqrt{ [l-m][l+m]} \; \ket{l-1,m},
\end{align}
with $\pi(a^*)$ being the hermitian conjugate of $\pi(a)$ and
$[x] := (q-q^{-1})^{-1}(q^x-q^{-x})$.
\end{proposition}
The proof is a long but straightforward calculation based on the
covariance property (\ref{equiv-1}) with the natural representation
$\rho$ of $\suq2$ on $\CH_h$ and the $\suq2$-module structure of
$\sq2$ given in (\ref{module}). The representations $\pi_\pm$ are
equivalent to the left regular representation of $\CA(SU_q(2))$ on
$L^2(SU_q(2))$ (with the Haar measure) when this representation is
restricted
to the subalgebra $\sq2$ and the representation space is restricted to the
($L^2$-completion) of certain vector spaces (left $\sq2$-modules)
\mbox{constructed in \cite{BM}.}

\bigskip
We take as the Hilbert space of our geometry $\CH = \CH_h \oplus \CH_h$,
with the natural grading $\gamma = \hbox{id} \oplus (-\hbox{id})$
and the representation $\pi(x) = \pi_+(x) \oplus \pi_-(x)$
for any $a \in \sq2$, which is equivariant with respect to the
diagonal action (which we call again $\rho$) of $\suq2$ on $\CH$.
Like in \cite{DaSi} we follow the method of equivariance to find first the real structure $J$. Let us
recall, that $J$ is the antiunitary part of an antilinear operator $T$
on $\CH$, which must then satisfy for any $h \in \suq2$ (on a dense
subspace of $\CH$), $\rho(h) T = T \rho(Sh)^*$.
Taking into account the required commutation relations with the
grading $\gamma$, that is $\gamma J = - J \gamma$ and $J^2=-1$, one
easily obtain that $J$ must be,
\begin{equation}
J \ket{l,m}_\pm = i^{2m} \ket{l,-m}_\mp \, , \label{defJ}
\end{equation}
where the label $\pm$ refers to the two copies of $\CH_h$ which
are marked by the eigenvalues of $\gamma$. With this data we
immediately meet an obstruction:
\begin{proposition}
The operator $J$ defined above does not satisfy the ``commutant"
requirement of a real spectral triple, that is, there
exist $x \in \sq2$ for which $J \pi(x) J^{-1}$ is not in
the commutant of $\pi(\sq2)$.
\end{proposition}
\smallskip
We move on to look for a variation of spectral geometry
up to infinitesimal. Let $\CK$ denotes
the ideal of compact operators on $\CH$ and $\CK_q \subset \CK$
be the ideal generated by operators $L_q$ of the form
$L_q \ket{l,m}_\pm = q^l \ket{l,m}_\pm$. Using compact perturbations
of the representation $\pi$ one proves the following 
\begin{proposition}
The operator $J$, defined in \eqref{defJ}, maps $\pi(\sq2)$ to its
commutant modulo compact operators (in fact modulo $\CK_q$). More
precisely, for any $x,y
\in \sq2$,
\begin{equation}\label{al-com}
\left[\pi(x), J \pi(y) J^{-1} \right]
\in \CK_q \, .
\end{equation}
\end{proposition}
Exact formul{\ae} shall be contained in the extended version of this note.

\bigskip
As a next step we derive the Dirac operator $D$. Beside postulating
that $D$ anticommutes with  $\gamma$ and commutes with $J$, we shall
also require that $D$ is equivariant, that is it commutes with the
representation $\rho$ of $\suq2$ on $\CH$. Each operator satisfying these
condition must be of the form
\begin{equation}
\label{dirac}
D \ket{l,m}_\pm =  d_l^\pm ~\ket{l,m}_\mp\,, \;\;\; d_l^\pm \in \R.
\end{equation}
\begin{proposition}
Up to rescaling and addition of an odd constant as well as of an odd 
element from $\CK_q$, there exists only one 
operator $D$ of the form (\ref{dirac})  which satisfies the order-one
condition up to compact operators (in fact modulo $\CK_q$), that is
for all
$x,y \in \sq2$,
\begin{equation}
\left[ J \pi(x) J^{-1}, [D, \pi(y)] \right] \in \CK_q\,.
\label{al-FOC}
\end{equation}
For this operator $D$, the parameters $d_l^\pm$ are given by
$d_l^+=d_l^- = l + \frac{1}{2}$.
\end{proposition}
\noindent
The condition (\ref{al-FOC}) has been verify explicitly for all
pairs of generators of $\sq2$ with the help of a symbolic
computation program. Furthermore, we have,
\begin{proposition}
For any $x \in \sq2$, the commutators $[D,\pi(x)]$ are bounded.
\end{proposition}
\medskip
It is evident that the operator $D$ is self-adjoint on a natural
domain in $\CH$ and that its resolvent is compact. Since
the spectrum of $|D|$ consists only of eigenvalues
$k=l+\frac{1}{2} \in{\N}$ with multiplicity $4k$, we managed
to realize the suggestion in \cite{CL} to use the $D$ with classical 
spectrum. Thus, the deformation
being isospectral, the dimension requirement is satisfied with
the spectral dimension of ($\sq2, \CH, D)$ being $n=2$.

We have made some advancement in the study of
the $q$-geometry and expect that similar structures exist on
other $q$-deformed spaces. There are still many points to be
addressed, notable the existence of a volume form and the
fulfillment of other axioms of spectral geometries. These points
shall be addressed in the extended version of the note.
An isospectral spectral triple for $SU_q(2)$ has been constructed in 
\cite{DLSSV}.

\end{document}